\documentclass[12pt,reqno]{article}
\usepackage{amsmath,amsthm,amsfonts,amssymb,amscd}
\usepackage[pdftex]{graphicx}
\usepackage{textcomp}
\usepackage{color}

\usepackage{pdfpages}

\usepackage[pdftex,bookmarks=true]{hyperref}
\usepackage[ansinew]{inputenc}
\usepackage{lscape}
\usepackage{multirow}
\usepackage{indentfirst} %(faz com que o 1ro paragrafo iniciando a secao seja igual aos demais ( no arquivo do Nassif dá erro)!
\usepackage{rotate}
\usepackage{marvosym}
\usepackage{eurosym}
\usepackage{psfrag}
\usepackage{color}

\theoremstyle{plain}
\newtheorem{thm}{Theorem}[section]

\newtheorem{clly}[thm]{Corollary}
\newtheorem{lemma}[thm]{Lemma}
\newtheorem{defi}[thm]{Definition}

\newtheorem{maintheorem}{Theorem}

\newcommand{\re}{{\Bbb R}}
\newcommand{\nat}{{\Bbb N}}

\title{Finiteness of homoclinic classes on sectional hyperbolic sets}
\author{A. M. L\'opez B, A.E. Arbieto
        }
\date{}

\begin{document}
\maketitle

\begin{abstract}

We study small perturbations of a sectional hyperbolic set 
of a vector field on a compact manifold. 
Indeed, we obtain robustly finiteness of  homoclinic classes 
on this scenary. Moreover, since attractor and repeller sets are particular cases 
of homoclinic classes, this result improve \cite{am}.

\end{abstract}

%%==================================================================SECTION 1======================================================================

\section{Introduction}

\noindent

A dynamical system looks for describe the future behavior on certain maps in different sceneries, 
namely, diffeomorphisms (discrete time) and vector fields (flows in continuous time). The properties 
of such future behavior, can be used for obtain global information as asymptotic behavior, stability among others. 
It is well known, that the concepts of sink and source, or attractor and repeller play an important role about the system's dynamic. 
Since these sets are a particular case of a homoclinic class, the study on this type of sets it becomes of more interest.

There are a lot of known results and examples it stress the importance of the sectional hyperbolic homoclinic 
classes, such as, the Spectral decomposition theorem and the classical construction of the geometric 
Lorenz models \cite{abs}, \cite{gw} among others. The study of sectional hyperbolic homoclinic classes is 
mostly open on the higher dimensional case.

The aim main of this paper is to study a problem related to the sectional hyperbolic sets for flows, that is,
finiteness of homoclinic classes in this scenary, and how many homoclinic classes can arise from small perturbations.

It's clear that there exist some motivations that come from the previous results. Specifically, in dimension three \cite{m}, \cite{ams}, 
show that the finiteness of attractors is obtained but for transitive or nonwandering flows, and this amount is in 
terms of the number of singularities. After, \cite{a} obtains the same results for transitive flows but in higher 
dimensions. Recently, \cite{am} improves these results proving finiteness of attractors and repellers for 
sectional hyperbolic sets in any dimension, and remove both transitivity and nonwandering hypotheses. 

However, some results talk about the Newhouse's phenomenon (coexistence of infinitely many sinks or sources) \cite{bodi}, 
or presence of infinitely homoclinic classes. In particular, \cite{carm} show that a generic $C^1$
vector field on a closed $n$-manifold has either infinitely many homoclinic classes, or a finite collection
of attractors (or, respectively, repellers) with basins that form an open-dense set.

Thus, it is natural to consider even the question if a sectional hyperbolic set 
has only a finite number of homoclinic classes. This question is related with some 
Palis and Bonatti's conjectures \cite{jp}, \cite{bo}. We provide a positive answer for this question.
In fact, we obtain robustly finiteness of  homoclinic classes 
on sectional hyperbolic sets. Moreover, since attractor and repeller sets are particular cases 
of homoclinic classes, this result improve \cite{am}.

Let us state our results in a more precise way.

Consider a compact manifold $M$ of dimension $n\geq 3$
with a Riemannian structure $\|\cdot\|$ (a {\em compact $n$-manifold} for short).
We denote by $\partial M$ the boundary of $M$.
Let ${\cal X}^1(M)$ be the space
of $C^1$ vector fields in $M$ endowed with the
$C^1$ topology.
Fix $X\in {\cal X}^1(M)$, inwardly
transverse to the boundary $\partial M$ and denote
by $X_t$ the flow of $X$, $t\in I\!\! R$.e flow and their extensions
The {\em maximal invariant} set of $X$ is defined by $M(X)= \displaystyle\bigcap_{t \geq 0} X_t(M).$
Notice that $M(X)=M$ in the boundaryless case $\partial M=\emptyset$.
A subset $\Lambda$ is called {\em invariant} if
$X_t(\Lambda)=\Lambda$ for every $t\in I\!\! R$.
We denote by $m(L)$ the minimum norm of a linear
operator $L$, i.e., $m(L)= inf_{v \neq 0} \frac{\left\|Lv\right\|}{\left\|v\right\|} $.

\begin{defi}
\label{d2}
A compact invariant set
$\Lambda$ of $X$
is {\em partially hyperbolic}
if there is a continuous invariant
splitting $T_\Lambda M=E^s\oplus E^c$ 
such that the following properties
hold for some positive constants $C,\lambda$:

\begin{enumerate}
\item
$E^s$ is {\em contracting}, i.e., 
$\mid\mid DX_t(x) \left|_{E^s_x}\right. \mid\mid
\leq Ce^{-\lambda t}$, 
for all $x\in \Lambda$ and $t>0$.
\item
$E^s$ {\em dominates} $E^c$, i.e., 
$\frac{\mid\mid DX_t(x) \left|_{E^s_x}\right. \mid\mid}{m(DX_t(x) \left|_{E^c_x}\right. )}
\leq Ce^{-\lambda t}$, 
for all $x\in \Lambda$ and $t>0$.
\end{enumerate}

We say the central subbundle $E^c_x$ of $\Lambda$ is 
{\em sectionally expanding} if
$$dim(E^c_x) \geq 2 \quad 
\mbox{and} \quad 
\left| det(DX_t(x) \left|_{L_x}\right. ) \right| \geq C^{-1}e^{\lambda t},
\mbox{ for all } x \in \Lambda,\quad t>0 $$ 
and all two-dimensional subspace $L_x$ of $E^c_x$.
Here $det(DX_t(x) \left|_{L_x}\right. )$ denotes
the jacobian of $DX_t(x)$ along $L_x$.
\end{defi}

Recall that a singularity of a vector field is hyperbolic if
the eigenvalues of its linear part
have non zero real part.

\begin{defi}
\label{shs}
A {\em sectional hyperbolic set}
is a partially hyperbolic set whose singularities are hyperbolic and whose central subbundle is sectionally expanding.
\end{defi}

The {\em $\omega$-limit set} of $p\in M$ is the set
$\omega_X(p)$ formed by those $q\in M$ such that $q=\lim_{n \rightarrow \infty}X_{t_n}(p)$ for some
sequence $t_n\to\infty$. The {\em $\alpha$-limit set} of $p\in M$ is the set
$\alpha_X(p)$ formed by those $q\in M$ such that $q=\lim_{n \rightarrow \infty}X_{t_n}(p)$ for some
sequence $t_n\to -\infty$. The {\em non-wandering set} of $X$ is
the set $\Omega(X)$ of points $p\in M$ such that
for every neighborhood $U$ of $p$ and every $T>0$ there
is $t>T$ such that $X_t(U)\cap U\neq\emptyset$.
Given $\Lambda \in M$ compact, we say that $\Lambda$ is {\em invariant} 
if $X_t(\Lambda)=\Lambda$ for all $t\in I\!\! R$.
We also say that $\Lambda$ is {\em transitive} if
$\Lambda=\omega_X(p)$ for some $p\in \Lambda$; {\em singular} if it
contains a singularity and
{\em attracting}
if $\Lambda=\cap_{t>0}X_t(U)$
for some compact neighborhood $U$ of it.
This neighborhood is often called
{\em isolating block}.
It is well known that the isolating block $U$ can be chosen to be
positively invariant, i.e., $X_t(U)\subset U$ for all
$t>0$.
An {\em attractor} is a transitive attracting set.
An attractor is {\em nontrivial} if it is
not a closed orbit.
A {\em repelling} is an attracting for the time reversed vector field $-X$ 
and a {\em repeller} is a transitive repelling set.

With these definitions we can state our main results.

\begin{maintheorem}
\label{thA}
For every sectional hyperbolic set $\Lambda$ of a vector field $X$ on a compact manifold there are
neighborhoods ${\cal U}$ of $X$, 
$U$ of $\Lambda$ and $n_0\in \mathbb{N}$ such that
\begin{center}
$\#\{L\subset U:$ $L$ is a Homoclinic class of $Y\in {\cal U}\} \leq n_0$.
\end{center}
\end{maintheorem}

% % % % % % % % % % % % % % % % % % % % % % % % % % % % % % % % % % % % % % % % % % % % % % % % % % % 
% % % % % % % % % % % % % % % % % % % % % % % % % % % % % % % % % % % % % % % % % % % % % % % % % % % 
% % % % % % % % % % % % % % % % % % % % % % % % % % % % % % % % % % % % % % % % % % % % % % % % % % % 
% % % % % % % % % % % % % % % % % % % % % % % % % % % % % % % % % % % % % % % % % % % % % % % % % % % 

Let us present some corollaries of Theorem \ref{thA}.
Recall that a {\em sectional Anosov flow} is a vector field whose maximal invariant set is sectional hyperbolic \cite{mem}.

\begin{clly}
\label{cor1}
For every sectional Anosov flow
of a compact manifold there are
a neighborhood $ {\cal U} \in {\cal X}^1(M)$ and $n_0 \in \mathbb{N}$ such that
\begin{center}
$\#\{L$ is a Homoclinic class of $Y\in {\cal U}\} \leq n_0$.
\end{center}
\end{clly}

%%============================================================================SECTION 3====================================================

\section{Preliminaries}
\label{sech}

\noindent
In this section, we recall some results on sectional hyperbolic sets, 
and we obtain some useful results for the main theorems. 

In first place, we recall the standard definition of hyperbolic set.

\begin{defi}
\label{hyperbolic}
A compact invariant set $\Lambda$ of $X$ is {\em
hyperbolic}
if there are a continuous tangent bundle
invariant decomposition
$T_{\Lambda}M=E^s\oplus E^X\oplus E^u$ and positive constants
$C,\lambda$ such that

\begin{itemize}
\item $E^X$ is the vector field's
direction over $\Lambda$.
\item $E^s$ is {\em contracting}, i.e.,
$
\mid\mid DX_t(x) \left|_{E^s_x}\right.\mid\mid
\leq Ce^{-\lambda t}$, 
for all $x \in \Lambda$ and $t>0$.
\item $E^u$ is {\em expanding}, i.e.,
$
\mid\mid DX_{-t}(x) \left|_{E^u_x}\right.\mid\mid
\leq Ce^{-\lambda t},
$
for all $x\in \Lambda$ and $t> 0$.
\end{itemize}
A closed orbit is hyperbolic if it is also hyperbolic, 
as a compact invariant set. An attractor is hyperbolic 
if it is also a hyperbolic set. 
\end{defi}

Let $\Lambda$ be a sectional hyperbolic set of a $C^1$ vector field $X$ of $M$.
Henceforth, we denote by $Sing(X)$ the set of singularities of the vector field $X$ and  
by $Cl(A)$ the closure of $A$, $A\subset M$. 

The following two results examining the sectional hyperbolic 
splitting $T_{\Lambda}M = E^s_{\Lambda} \oplus E^c_{\Lambda}$
of a sectional hyperbolic set $\Lambda$ of $X \in {\cal X}^1(M)$
appear in \cite{mpp2} %%\cite[Lemma 3,Pag 5]{mpp2}, \cite[Theorem A,Pag 3]{mpp2}
for the three-dimensional case and in \cite{lec} %%\cite[Corollary 2.7, Pag 65]{lec} , \cite[Lemma 2.7, Pag 67]{lec}
for the higher dimensional case.

\begin{lemma}
\label{l1}
Let $\Lambda$ be a sectional hyperbolic set of a $C^1$ vector field $X$ of $M$.
Then, there are a neighborhood $ {\cal U} \subset {\cal X}^1(M)$ of $X$ 
and a neighborhood $U \subset M$ of $\Lambda$ such that if $Y \in {\cal U}$,
every nonempty, compact, non singular, invariant set $H$
of $Y$ in $U$ is hyperbolic {\em saddle-type} (i.e. $E^s\neq 0$ and $E^u\neq 0$).
\end{lemma}

\begin{thm}
\label{t1}
Let $\Lambda$ be a sectional hyperbolic set of a $C^1$ vector field $X$ of $M$.
If $\sigma\in Sing(X)\cap \Lambda$, then
$\Lambda \cap W^{ss}_X(\sigma)=\{\sigma\}.$
\end{thm}

Next we explain briefly how to obtain sectional hyperbolic sets 
nearby $\Lambda$ from vector fields close to $X$. 
Fix a neighborhood $U$ with compact closure of $\Lambda$ as in Lemma \ref{l1}.
Define
$$\Lambda_X = \cap_{t \in \re} X_t(Cl(U)).$$ 
Note that $\Lambda_X$ is sectional hyperbolic and $\Lambda \subset \Lambda_X$.
Likewise, if $Y$ is a $C^1$ vector field close to $X$, we define the continuation
$$\Lambda_{Y} = \cap_{t \in \re} Y_t(Cl(U)).$$

%%%%%%%%%%%%%%%%%%%%%%%%%%%%%%%%%%%%%%%%%%%%%%%%%%%%%%%%%%%%%%%%%%%%%%%%%%%%%%%%%%%%%%%%%%%%%%%%%%%%%%%%%%%
%%%%%%%%%%%%%%%%%%%%%%%%%%%%%%%%%%%%%%%%%%%%%%%%%%%%%%%%%%%%%%%%%%%%%%%%%%%%%%%%%%%%%%%%%%%%%%%%%%%%%%%%%%%
%%%%%%%%%%%%%%%%%%%%%%%%%%%%%%%%%%%%%%%%%%%%%%%%%%%%%%%%%%%%%%%%%%%%%%%%%%%%%%%%%%%%%%%%%%%%%%%%%%%%%%%%%%%

\subsection{Linear Poincaré flow and scaled Poincaré flow}

Recall, ${\cal X}^1(M)$ is the space
of $C^1$ vector fields in $M$ endowed with the
$C^1$ topology, $X\in {\cal X}^1(M)$ a fixed vector field, inwardly
transverse to the boundary $\partial M$. We denote
by $X_t: M \rightarrow M$ the flow generated of $X$, $t\in I\!\! R$, and the tangent flow 
$DX_t : TM \rightarrow TM$. 

Denote the normal bundle of $X$ by $$N=N^X= \bigcup \limits_{x \in M^{*}} N_x$$ 
where $N_x$ represents the orthogonal complement of the flow direction $X(x)$, i.e.,
$N_x = \{v \in T_xM : v \perp X(x)\}$ and $M^*= M \setminus Sing(X)$.

Thus, given $x \in M^*$ and $v \in  N_x$, we denote by $P_t(x)v$ the orthogonal projection of $DX_(x)(v)$ on $N_{X_t(x)}$ along
the flow direction, i.e.,
$$P_t(x)v = DX_t(x)v - \frac{\langle DX_t(x)(v), X(X_t(x))\rangle}{||X(X_t(x))||^2} X(X_t(x)),$$
where $\langle \cdot, \cdot \rangle $ is the inner product on $T_xM$ given by the Riemannian metric.

The flow $P_t$ is called the "linear Poincaré flow (LPF)`` and exhibits certain properties 
such as $||P_t||$ is uniformly bounded for $t$ in any bounded interval although it is just
defined on the regular set which is not compact in general. 
% % This flow could also be defined in a more general way by Liao used the terminology
% % of ?extended linear Poincar ?e flow?.
In the same way of \cite{ls}, \cite{gy}, we shall use another useful flow generated by the LPF. 
Thus, given $x\in M^*$, $v \in N_x$ we denote by  $P_t^*: N \rightarrow N$ (called scaled linear Poincaré flow) 
and defined by 
$$P_t^*(x)v= \frac{\|X(x)\|}{\|X(X_t(x))\|} P_t(x)v= \frac{P_t(x)v}{\|DX_t(x)|_{\langle X(x) \rangle} \|}$$

Here, this scaled linear Poincaré flow will allow us to work with some difficulties produced
by singularities.

The next lemma states the basic properties of star flows, proved in \cite{ls1}.

\begin{lemma}\label{mixto}%[WE CAN CHOOSE BY STAR FLOW]
For any sectional Anosov flow $X$, $X \in \mathcal{X}^1(M)$, there are a $C^1$-neighborhood $U$ and
numbers $\eta > 0$ and $T > 0$ such that for any periodic orbit $O_Y$ of $Y \in U$ with period
$\pi(O_Y) \geq T$, if $N_{O_Y} = N^s \oplus N^u$ is the hyperbolic splitting with respect to $P_t^Y$ then

\begin{itemize}
 \item For every $x \in O_Y$ and $t \geq T$, one has
 \begin{equation}
    \begin{tabular}{l}\\ 
      \large{ $\frac{\|P_t^Y |_{N^s(x)} \|}{m(P_t^Y |_{N^u(x)})} \leq e^{-2\eta t}$}\\
    \end{tabular}
  \label{mix1}
 \end{equation}

\item If $x\in O_Y$, then
\begin{equation}
  \begin{tabular}{l}\\
  $ \prod \limits_{i=0} ^{[ \frac{\pi(O_Y)}{T}] -1} \|P_T^Y |_{N^s(Y_it(x))} \| \leq e^{-\eta \pi(O_Y)}$\\
  $ \prod \limits_{i=0} ^{[ \frac{\pi(O_Y)}{T}] -1} m (P_T^Y |_{N^u(Y_it(x))} ) \geq e^{\eta \pi(O_Y)}$ \\
  \end{tabular}
\label{mix2}
\end{equation}

\end{itemize}

\end{lemma}

Here $m(A)$ is the mini-norm of $A$, i.e., $m(A)= \|A^{-1}\|^{-1}$.

%%It gives uniform estimations on some non-compact sets.

\begin{defi}
\label{etacont}
Let $X \in  \mathcal{X}^1 (M)$, $\Lambda$ a compact invariant set of $X$ , 
and $E \subset N_{\Lambda \setminus Sing(X)}$ an invariant bundle of the linear Poincaré flow $P_t$. For $\eta > 0$ and
$T > 0$, $x \in \Lambda \setminus Sing(X)$ is said to be $(\eta, T, E )^*$-contracting if for any $n \in \nat$, 
$$ \prod \limits_{i=0} ^{n-1} \|P_T^* |_{E(X_it(x))} \| \leq e^{-\eta n }$$

Similarly, $x \in \Lambda \setminus Sing(X)$ is said to 
be $(\eta, T, E)^*$-expanding if it is $(\eta, T, E )^*$-contracting for $-X$.

\end{defi}

The following theorem provides existence of local invariant stable and unstable manifolds, 
with size proportional to the velocity of the vector field.

\begin{thm}\label{e}

Let $X \in \mathcal{X}^1(M)$ and $\Lambda$ a compact invariant set of $X$ . Given
$\eta > 0$, $T > 0$, assume that $N_{\Lambda \setminus Sing(X)} = E \oplus F$ is an $(\eta, T )$-dominated splitting with
respect to the linear Poincaré flow. Then, for any $\epsilon > 0$, there is $\delta > 0$ such that if
$x$ is $(\eta, T, E )^*$ contracting, then there is a $C^1$ map $k : E_x (\delta \|X(x)\|) \rightarrow N_x$ such that

\begin{itemize}
 \item $d_{C^1}(k,Id) < \epsilon$
 
 \item $k(0) =0$
 
 \item $W^{cs}_{\delta \| X(x)\|}(x) \subset W^s(O(x))$, where $W^{cs}_{\delta \| X(x)\|} = exp_x (Im(k))$ 
 
 %\item $W^{cs}_{\delta \| X(x)\| \subset W^s(O(x))$, where $W^{cs}_{\delta \| X(x)\| = exp_x (Im(k))$ 
\end{itemize}

Here $E_x(r) = {v \in E_x : \|v\| \leq r}.$
 
\end{thm}

It is well known that is possible found an uniform size of stable manifolds 
for diffeomorphisms and nonsingular flows exhibiting an uniform contraction on its derivate. 
Thus, the above result for the singular case, provides a construction of 
size of stable manifolds being proportional to the flow speed.

% % % % % % % % % % % % % % % % % % % % % % % % % % % % % % % % % % % % % % % % % % % % % % % % % % % 
% % % % % % % % % % % % % % % % % % % % % % % % % % % % % % % % % % % % % % % % % % % % % % % % % % % 
% % % % % % % % % % % % % % % % % % % % % % % % % % % % % % % % % % % % % % % % % % % % % % % % % % % 

%%%===============================================================================================

\subsection{Lorenz-like singularity and singular cross-section}

Let $M$ be a compact $n$-manifold, $n \geq 3$.
Fix $X\in {\cal X}^1(M)$, inwardly 
transverse to the boundary $\partial M$. We denote
by $X_t$ the flow of $X$, $t\in I\!\! R$.\\

By using the standard definition of hyperbolic set (Definition \ref{hyperbolic}),
it follows from the stable manifold theory \cite{hps} that if $p$ belongs to a hyperbolic set $\Lambda$, then the following sets

\begin{tabular}{lll}
$W^{ss}_X(p)$ & = & $\{x:d(X_t(x),X_t(p))\to 0, t\to \infty\}$ and\\
$W^{uu}_X(p)$ & = & $\{x:d(X_t(x),X_t(p))\to 0, t\to -\infty\}$ \\
\end{tabular}\\
are $C^1$ immersed submanifolds of $M$, which are tangent at $p$ to the subspaces $E^s_p$ and $E^u_p$ of $T_pM$, respectively.
Similarly,

\begin{tabular}{lll}
$W^{s}_X(p)$ & = & $ \bigcup_{t\in I\!\! R}W^{ss}_X(X_t(p))$ and\\
$W^{u}_X(p)$ & = & $ \bigcup_{t\in I\!\! R}W^{uu}_X(X_t(p))$ \\
\end{tabular}\\
are also $C^1$ immersed submanifolds tangent to $E^s_p\oplus E^X_p$ and $E^X_p\oplus E^u_p$ at $p$, respectively.
Moreover, for every $\epsilon>0$ we have that

\begin{tabular}{lll}
$W^{ss}_X(p,\epsilon)$ & = & $\{x:d(X_t(x),X_t(p))\leq\epsilon, \forall t\geq 0\}$ and\\
$W^{uu}_X(p,\epsilon)$ & = & $\{x:d(X_t(x),X_t(p))\leq \epsilon, \forall t\leq 0\}$\\
\end{tabular}\\
are closed neighborhoods of $p$ in $W^{ss}_X(p)$ and $W^{uu}_X(p)$, respectively.

For simplicity, given $\epsilon > 0$,  the above submanifolds will be denoted by 
$W^{ss}_{loc}(p)$, $W^{uu}_{loc}(p)$, $W^{s}_{loc}(p)$ and $W^{uu}_{loc}(p)$ 
respectively.\\

It is well known from stability theory for hyperbolic sets, that we can fix a neighborhood 
$U \subset M$ of $\Lambda$, a neighborhood ${\cal U} \subset {\cal X}^1(M)$ 
of $X$ and $\epsilon>0$ such that every hyperbolic set $H$ in $U$ of 
every vector field $Y$ in $\mathcal{U}$ satisfies that its local stable and instable manifold  
\begin{equation}
	\begin{tabular}{l}
		$W^{ss}_Y(x,\epsilon)$ and $W^{uu}_Y(x,\epsilon)$ have uniform size $\epsilon$
		for all $x \in H.$\\
	\end{tabular}
\label{El2}
\end{equation}

There is also a stable manifold theorem in the case when $\Lambda$ is a sectional hyperbolic set.
Indeed, if we denote by $T_{\Lambda}M=E^s_{\Lambda}\oplus E^c_{\Lambda}$ the corresponding sectional hyperbolic
splitting over $\Lambda$, we assert that there exists such contracting foliation on a small neighborhood $U$ 
of $\Lambda$. Note that this extended foliation is not necessarily invariant, and we can only ensure the 
invariance if this one, is at least, an attracting set \cite{appv}. %%\cite[Section 2.1]{appv}
This extension is carried out as follows: first, we can 
choose cone fields on $U$ and we consider the space of tangent foliations to the cone fields. 
Given a point $x \in U$, whenever the positive orbit remain within to $U$, for example 
$t=1$, we can use the derivate map $DX_{-1}(x)$. This map sends the leaf at $X_{-1}(x)$ 
inside of cone $C_{X_{-1}(x)}$ to the leaf at $x$ inside of cone $C_x$, 
contracting the angle and stretching the tangent vectors to the initial foliation. Then, we can apply 
fiber contraction \cite{hps}, \cite{chi}. %%\cite[Pag 30,31,80]{hps}, \cite[Theorem 1.243, Pag 127]{chi}
Now, by using the Fiber Contraction Theorem \cite{hps} 
the foliation arises. Thus, we have from \cite{hps} that
the contracting subbundle $E^s_{\Lambda}$
can be extended to a contracting subbundle $E^s_U$ in $M$ (not necessarily invariant). 

Moreover, such an extension by construction is tangent to a continuous foliation denoted 
by $W^{ss}$ (or $W^{ss}_X$ to indicate dependence on $X$).
By adding the flow direction to $W^{ss}$ we obtain a continuous foliation 
$W^s$ (or $W^s_X$) now tangent to $E^s_U\oplus E^X_U$.
Unlike the hyperbolic case $W^s$ may have singularities, all of which being
the leaves $W^{ss}(\sigma)$ passing through the singularities $\sigma$ of $X$.

It turns out that every singularity $\sigma$ of a sectional hyperbolic set 
$\Lambda$ satisfies $W^{ss}_X(\sigma)\subset W^s_X(\sigma)$.
Furthermore, there are two possibilities for such a singularity, namely,
either $dim(W^{ss}_X(\sigma))=dim(W^s_X(\sigma))$ 
(and so $W^{ss}_X(\sigma)=W^s_X(\sigma)$) or $dim(W^{s}_X(\sigma))=dim(W^{ss}_X(\sigma))+1$. 
In the later case we call it Lorenz-like according to the following definition. 

\begin{defi}
\label{ll} 
Let $\Lambda$ be a sectional hyperbolic set of a $C^1$ vector field $X$ of $M$.
We say that a singularity $\sigma$ of $\Lambda$ is {\em Lorenz-like}
if
$dim (W^s(\sigma))=dim (W^{ss}(\sigma))+1.$
\end{defi}

Hereafter, we will denote $dim(W^{ss}_X(\sigma))=s$, $dim(W^{u}_X(\sigma))=u$ and 
therefore $dim(W^{s}_X(\sigma))=s+1$ by definition. 
Moreover, $W^{ss}_X(\sigma)$ separates $W^s_{loc}(\sigma)$
in two connected components which we will denote by $W^{s,t}_{loc}(\sigma)$ 
and $W^{s,b}_{loc}(\sigma)$, respectively 
(in this way $W^{s,t}_{\cdot}(\cdot)$ and $W^{s,b}_{\cdot}(\cdot)$ denote top and bottom stable components).

\section{Finiteness}

We start by recalling some useful definitions to prove 
the lemmas and the propositions that 
provide very important properties on sectional hyperbolic sets, 
that in our case support the main theorems' proofs.

Let $O = \left\{X_t (x): t \in \re \right\}$ be the orbit of $X$ through $x$,
then the stable and unstable manifolds of $O$ defined by
\begin{center}
	$W^s(O) =\cup_{x\in O} W^{ss}(x)$, and $W^u(O) =\cup_{x\in O} W^{uu}(x)$
\end{center}
are $C^1$ submanifolds tangent to the subbundles $E^s_{\Lambda} \oplus E^X_{\Lambda}$ 
and $E^X_{\Lambda} \oplus E^u_{\Lambda}$, respectively.

A {\em homoclinic orbit} of a hyperbolic periodic orbit $O$ is an orbit in $\gamma \subset W^s(O) \cap W^u(O)$.
If additionally $T_qM = T_qW^s(O) + T_qW^u(O)$ for some (and hence all) point $q \in \gamma$, then we say that
$\gamma$ is a {\em transverse homoclinic orbit} of $O$.  

\begin{defi}
The {\em homoclinic class} $H(O)$ of a hyperbolic periodic orbit $O$ is
the closure of the union of the transverse homoclinic orbits of $O$.
We say that an invariant set $L$ is a {\em homoclinic class} 
if $L = H(O)$ for some hyperbolic periodic orbit $O$.
\end{defi}

Recall that $Sing(X)$ denotes the set of singularities of the vector field $X$ and 
$Cl(A)$ denotes the closure of $A$, $A\subset M$. Moreover, 
for $\delta>0$ we define $B_\delta(A)=\{x\in M:d(x,A)<\delta\}$, 
where $d(\cdot,\cdot)$ is the metric in $M$. \\

The following two results examining the finiteness of homoclinic classes on a sectional hyperbolic set in 
$\Lambda$ of $X \in {\cal X}^1(M)$, in certain particular cases appear in \cite{am}.\\

\begin{thm}
\label{fini}
For every sectional hyperbolic set $\Lambda$ of a vector field $X$ on a compact manifold there are
neighborhoods ${\cal U}$ of $X$, 
$U$ of $\Lambda$ and $n_0\in \mathbb{N}$ such that
\begin{center}
$\#\{L\subset U:$ $L$ is an attractor or repeller of $Y\in {\cal U}\} \leq n_0$.
\end{center}
\end{thm}

\bigskip

\begin{lemma}
\label{le2}
Let $X$ be a $C^1$ vector field
of a compact $n$-manifold $M$, $X \in {\cal X}^1(M)$. 
Let $\Lambda \in M$ be a hyperbolic set of $X$.
Then, there are a neighborhood $ {\cal U} \subset {\cal X}^1(M)$ of $X$, 
a neighborhood $U \subset M$ of $\Lambda$ and $n_0 \in \mathbb{N} $ such that
\begin{center}
	 $\#\{L\subset U:$ $L$ is homoclinic class of $Y\in {\cal U}\} \leq n_0$
\end{center}
for every vector field $Y \in {\cal U}$.
\end{lemma}

%%=======================================================================SECTION 4============================================================

\section{Proof of the main theorems}

\subsection{Proof of Theorem \ref{thA}}

%%\textbf{Proof of Theorem \ref{thA}}
\begin{proof}

We prove the theorem by contradiction.
Let $X$ be a $C^1$ vector field
of a compact $n$-manifold $M$, $n\geq3$, $X \in {\cal X}^1(M)$. 
Let $\Lambda \in M$ be a sectional hyperbolic set of $X$.
Then, we suppose that
there is a sequence of vector fields $(X^n)_{n\in \nat} \subset {\cal X}^1(M)$, 
$X^n\overset{C^1}{\to} X$ such that every vector field $X^n$ exhibits 
$n$ homoclinic classes, with $n>n_0$. 

It follows from the theorem \ref{fini} 
that there are neighborhoods $ {\cal U} \subset {\cal X}^1(M)$ of $X$ 
and $U \subset M$ of $\Lambda$
such that the attractor and repellers in $U$ are finite for all vector field $Y$ in ${\cal U}$.
Then, at least most of them are of saddle type and with loss of generality we can assume all of saddle type. 
Thus, we are left to prove only for this case. 

We denote by $L^n$ an homoclinic class of $X^n$ in $\Lambda_{X^n}$. 
Since $\Lambda_{X^n}$ is arbitrarily close to $\Lambda_X$ and $L^n \in \Lambda_{X^n}$, 
$L^n$ is also arbitrarily close to $\Lambda_X$. Therefore, 
we can assume that $L^n$ belongs to $\Lambda_X$ for all $n$.

Let $(L^n)_{n\in \nat}$ be the sequence of  
homoclinic classes contained in $\Lambda_X$. We assert that  
\begin{center}
$Sing(X) \bigcap (\cap_{N >0} Cl(\cup_{m\geq N} L^n)) \neq \emptyset. $
\end{center}

Otherwise, as in \cite{am}, we can find $\delta>0$, such that 
$$B_\delta(Sing(X))\bigcap\left(\cup_{n\in \nat} L^n\right)=\emptyset.$$

Thus, we define
$H=\bigcap_{t\in I\!\! R}X_t\left(U\setminus B_{\delta/2}(Sing(X))\right)$.
%%By definition $Sing(X)\cap H=\emptyset$, $H$ is compact since $\Lambda$ is, 
%%and $H$ is a nonempty compact set \cite{a}, which is clearly invariant for $X$.
It follows from the Lemma \ref{l1} that $H$ is a hyperbolic set 
and beside the Lemma \ref{le2} we have that 
there are neighborhoods $ {\cal U} \subset {\cal X}^1(M)$ of $X$, 
$U \subset M$ of $H$ and $n_1 \in \mathbb{N} $ such that
\begin{center}
	 $\#\{L\subset U:$ $L$ is an homoclinic class of $Y\in {\cal U}\} \leq n_0.$
\end{center}
Note that the last inequality holds for every vector field $Y \in {\cal U}$, 
but this leads to a contradiction, since by hypothesis 
we have that
\begin{center}
$\#\{L\subset H:$ $L$ is an homoclinic class of $Y\in {\cal U}\} \geq n>n_0$.
\end{center}

Thus, we have that there exists $p_n \in L_n$ periodic point for the vector field $X^n$ such that 
$p_n \rightarrow \sigma$, with $\sigma \in \Lambda$ a singularity of $X$. 
Indeed, for each vector field $X^n$ we have at least a sequence of $n$ periodic points $p_1^n, p_2^n,..., p_n^n$ 
that represents each homoclinic class in $\Lambda_X$. We shall write $p_n$ for short and we will clarify if necessary.

From the Lemma \ref{mixto} there exist neighborhood $U$ and numbers $\eta, T >0$ such that each point $p_n$ 
satisfies the inequalities (\ref{mix1}) and (\ref{mix2}) \cite{sgw}, \cite{ls}. 

By Poincaré recurrence there is $x_n \ R(X_T)$ satisfying the following inequality.

\begin{equation}
\label{teles}
 \lim \limits_{n \rightarrow \infty} \frac{1}{nT} (log X(X _{nT}(x_n))- log X(x_n)) = 0
\end{equation}

Then, it follows from \ref{teles} that the orbit $O_{X_n}=O_{X_n}(x_n)= O_{X_n}(p_n)=O_n$ 
through $p_n$ is eventually $(\eta, T)^*$-contractible with reference point $x_n$, with respect to scaled LPF.

Thus, we have a sequence of periodic points $(p_n)_{n \in \nat}$ such that:
\begin{itemize}
 \item $p_n \to \sigma.$
 \item $p_i$ is not homoclinic related with $p_j$ for every $i,j \in \nat.$
 \item $x_n$ is $(\eta, T)^*$-contractible with respect to scaled LPF, for all $n \in \nat.$
\end{itemize}

Without loss of generality, we can suppose that the sequence $(x_n)_{n\in \nat}$ is far from of set $Sing(X)$. 
Since $p_n \to \sigma$, passing to a subsequence if necessary, 
we can suppose that this one sequence converges to point $z \in W^s_X(\sigma)$.

It follows of the definition \ref{etacont} and Theorem \ref{e}, that each $x_n$ is endowed with a 
unstable manifold whose size is proportional to velocity of vector field, i.e., given $\epsilon >0$ there is $\delta >0$ 
such that for each $x_n$ we have that 
there is a 
$$W^{cu}_{\delta \| X(x_n)\|}(x_n) \subset W^u(O(x_n)).$$%% where $W^{cu}_{\delta \| X(x)\|} = exp_x (Im(k))$

By continuous parameters variation, one has that 

$$A^n = \bigcup_{0\leq t \leq 1} X_t(W^{cu}_{\delta \| X(x)\|}(x_n)) \subset W^u(O(x_n)).$$

Since $x_n \to z$, for $n$ enough large, there is $N \in \nat$ such that if $i,j \geq N$ one has 
$$d(x_i, x_j) < \frac{\delta}{\|X(z)\|} \quad and \quad X(z) \approx X(x_i) \approx X(x_j)$$

Thus, we have proved that $W^u(O(x_n))$ contains a local unstable manifold proportional to the flow speed. 
On the other hand, it implies that $$A^i \cap A^j \neq \emptyset,$$ 
but it implies that $p_i$ is homoclinic related with $p_j$ too. This is a contradiction and the proof follows.

\end{proof}

\hspace{1.0 cm}

%%============================================================================SECTION 3====================================================

\subsection{Proof of the corollaries}

\textbf{Proof of Corollary \ref{cor1}}

\begin{proof}
Since $X$ is a sectional Anosov flow, then its maximal invariant $M(X)$ is 
a sectional hyperbolic set for $X$. By using the Theorem \ref{thA} for $M(X)$ 
the proof follows.

\end{proof}

%%================================================================================================================00
%%=================================================THE ARTICLE IS FINISH =============================================================00
%%=================================================================================================================00

% % % \bibliographystyle{ieeetr}
% % % \bibliography{RCMBibTeX}

\bibliographystyle{acm}

% \bibliography{RCMBibTeX}

\medskip 

\flushleft
A. M. L\'opez B\\
Instituto de Matem\'atica, Universidade Federal Rural do Rio de Janeiro\\
Rio de Janeiro, Brazil\\
E-mail: andresmlopezb@gmail.com

A. E. Arbieto \\
Instituto de Matem\'atica, Universidade Federal do Rio de Janeiro\\
Rio de Janeiro, Brazil\\
E-mail: arbieto@im.ufrj.br

\end{document}